\documentclass[12pt]{article} 
\usepackage{graphicx,amsmath, amssymb, verbatim}
\usepackage[all]{xy}
\newtheorem{proposition}{Proposition}

\newtheorem{conjecture}{Conjecture}

\newcommand {\C } {\mathbb{C}} 

\newcommand {\p } {\mathbb{P}} 
 
\newcommand {\Z} {\mathbb{Z}}
\newcommand {\oo} {\mathcal O}
\newcommand {\la} {\lambda}

\newcommand{\ba}{\begin{eqnarray}}
\newcommand{\ea}{\end{eqnarray}}
\newcommand{\no}{\nonumber}

\begin{document} 

\title{\bf{Local mirror symmetry of curves: \\ Yukawa couplings
and genus 1}}

\author{Brian Forbes \\
\it Research Institute for Mathematical Sciences \\ \it  Kyoto University \\ \it  Kyoto 606-8502, Japan \\ \it{brian@kurims.kyoto-u.ac.jp}\\ \\
Masao Jinzenji  \\ \it Division of 
Mathematics, Graduate School of Science \\ \it Hokkaido University \\ 
\it  Sapporo  060-0810, Japan\\
\it jin@math.sci.hokudai.ac.jp }

\maketitle 

\begin{abstract}
We continue our study of equivariant local mirror symmetry of curves,
 i.e. mirror symmetry for $X_k=\oo(k)\oplus \oo(-2-k)\rightarrow \p^1$
 with torus action $(\la_1,\la_2)$ on the bundle. For the
 antidiagonal action $\la_1=-\la_2$, we find closed formulas for the
 mirror map, a rational $B$ model Yukawa coupling and consequently
 Picard-Fuchs equations for all
 $k$. Moreover, we give a simple closed form for the $B$ model genus 1
 Gromov-Witten potential. For the diagonal action $\la_1=\la_2$, we
 argue that the mirror symmetry computation is equivalent to that of the
 projective bundle $\p(\oo\oplus \oo(k)\oplus \oo(-2-k))\rightarrow \p^1$. Finally, we outline the computation of  equivariant Gromov-Witten invariants for $A_n$ singularities and toric tree examples via mirror symmetry.
\end{abstract}

\section{Introduction}

Mirror symmetry has for some time
now provided a convenient shortcut in  the computation of Gromov-Witten
invariants of toric varieties. Although mirror symmetry for compact
Calabi-Yau toric varieties historically appeared first, it has since been realized
that the noncompact (or local) CY toric variety setting enjoys many
simplifications not present in the compact case. 

However, as it turns out, there is a price to be paid for these
simplifications: when working with noncompact Calabi-Yaus, one loses the
nice structure theorems regarding periods of the mirror Calabi-Yau
manifold. For example, on compact threefolds, the constant term of
the prepotential (generating function of genus 0 Gromov-Witten
invariants) is known to be given by the triple intersection number of
curve classes in the space. However, for noncompact spaces, we have no
triple intersection number, and often the natural choice for this number
(from the perspective of mirror symmetry)
turns out to be fractional.   

More along the lines of the present paper, we are unable to use mirror
symmetry for any Calabi-Yau threefold $X$ such that $b_4(X)=0$. This is
because the periods mirror to the four cycles of $X$ are used to compute
the prepotential. As an even more serious problem, consider the case of `local
mirror symmetry of $\p^1$', by which we mean mirror symmetry for
$X_k=\oo(k)\oplus\oo(-2-k)\rightarrow \p^1$. The reason for this
terminology is that if $M$ is a Calabi-Yau threefold containing an
imbedded $\p^1$, $\p^1\hookrightarrow M$, then the local Calabi-Yau
condition implies
\ba
N_{\p^1/X}\cong \oo(k)\oplus\oo(-2-k).
\ea
This space is actually greatly troubled, as the Gromov-Witten theory of
$X_k$ for $k\ge 0$ is not even well-defined! This was shown in
\cite{BP}, where it was found that the Gromov-Witten invariants change
drastically with different choices of equivariant weights
$(\la_1,\la_2)$ acting on the bundle. As such, any mirror symmetry construction for
$X_k$ will necessarily be one of \it equivariant \normalfont mirror symmetry.

Thus, in this paper we develop a version of mirror symmetry for $X_k$
with torus weights $(\la,-\la)$ acting on the bundle. This case is of special
interest, as this is the Gromov-Witten theory computed by physicists
\cite{M}. We show that through a certain decomposition of the bundle, we
are able to describe mirror symmetry at genus 0 (i.e., the mirror map
and Yukawa coupling) via very simple
rational functions. Moreover, this structure allows us to easily work out mirror
symmetry at genus 1.

After this paper was completed, we were informed by the authors of
\cite{M} that they had obtained the same formula for the mirror map and
genus 1 partition function in their paper. In contrast to our work here, the computations of
\cite{M} are ultimately from the topological vertex, or $A$ model,
perspective, and their mirror map was found as a `natural' variable for
the problem, rather than a mirror map in the strict sense. Here, we will perform all calculations using only
the techniques of
mirror symmetry- and we do find the same mirror map can be used, though this mirror
map actually belongs to a different space, as we will see.

The organization of this paper is as follows. Section 2 summarizes our
main results. In Section 3, we review our previous results \cite{FJ3},
and give a derivation of natural rational Yukawa couplings
and genus 1 mirror symmetry
on $X_k$ (with antidiagonal action $(\la,-\la)$). Finally, in Section 4
we apply equivariant mirror symmetry to two spaces which are not of the
bundle-over-$\p^1$ type.

\bigskip

\bigskip

\bf{Acknowledgements}\normalfont

The authors would like to thank A. Givental, H. Iritani, S. Hosono,
M. Marino, K. Saito and A. Takahashi for helpful discussions. We also thank
Y. Konishi for providing us with a computer program for the topological
vertex calculation. We also thank N. Caporaso, L. Griguolo, M. Marino,
S. Pasquetti and D. Seminara for drawing our attention to the content of \cite{M}.

\section{Overview}

We state our main results and methods in this section. Throughout, we take $X_k=\oo(k)\oplus \oo(-2-k)\rightarrow \p^1$ and equip the bundle with a torus action with weights $(\la_1,\la_2)$. Let $t_k$ to be the K\"ahler class of $\p^1\hookrightarrow X_k$.

\subsection{Curves.}

We consider first the antidiagonal action $\la_1=-\la_2$ case. This is the computation which is of interest to physicists, and is the so-called equivariantly Calabi-Yau setting. We can exhibit the Calabi-Yau property by observing that the sum of the column vectors of the matrix of charge vectors for $X_k$ is zero when $\la_1=-\la_2$:
\ba
\label{matrix0}
\begin{pmatrix}
1&1&k&-2-k \\
0&0&\la_1&\la_2
\end{pmatrix}
\ea
Mirror symmetry was of course first observed for Calabi-Yau manifolds, and many of the nice structures associated to quantum cohomology, etc. owe themselves to this property. As such, although equivariant Gromov-Witten invariants for $X_k$ have been to this point not well understood, one may hope that equivariantly Calabi-Yau spaces possess the same structure. We have found that in fact this is the case:  
\begin{conjecture}
For  $X_k=\oo(k)\oplus \oo(-2-k)\rightarrow \p^1$ with antidiagonal torus action
 $(\la,-\la)$ on the bundle, the mirror map is given by 
\ba
q\frac{dt_k}{dq}=\frac{1+(-1)^{k+1}(k+1)^2 q}{1+(-1)^{k+1}q}
\ea
and the rational $B$ model Yukawa coupling, by the formula
\ba
Y^k_{qqq}=\frac{-1}{k(k+2)}\Big(q\frac{dt_k}{dq}\Big)^2.
\ea
In particular, this implies that we have the following Picard-Fuchs
 equation describing mirror symmetry for $X_k$:
\ba
\partial_{t_1}\Big(\frac{1}{Y^k_{t_1t_1t_1}}\Big)\partial_{t_1}^2
 =\theta^2\Big(q\frac{dt_k}{dq}\Big)^{-1}\theta.
\ea
\end{conjecture}
As mentioned in the introduction, the mirror map above was found in
\cite{M} as a `natural' variable for this calculation, and it was then
speculated in \cite{M} based on integrality properties that this formula
might be interpreted as a mirror map. Here, we have found directly from
mirror symmetry that this is indeed the mirror map, though technically, the
mirror map of a different space (see Proposition 1 below).  

We note in particular that this formula implies that the constant term of the $A$ model Yukawa coupling will be fractional:
\ba
Y^k_{t_kt_kt_k}=\frac{\partial^3 \mathcal F_k}{\partial t_k^3}=\Big(q\frac{dt_k}{dq}\Big)^{-3}Y^k_{qqq}=\frac{-1}{k(k+2)}+O(e^t).
\ea
where $\mathcal F_k$ is the prepotential for $X_k$.
This is precisely the value that was predicted (through entirely
different considerations) in \cite{AOSV}.
 As this constant has
the interpretation of being the triple intersection number of
$\p^1\hookrightarrow X_k$, we see the fractionality of intersection
numbers which was observed previously in \cite{JN},\cite{FJ1},
\cite{FJ2}. We note that this is choice of triple intersection number is
not the unique one that gives a rational Yukawa coupling; however, this choice
gives the simplest form of the the $B$ model Yukawa
coupling and Picard-Fuchs equations, hence making it natural from the
$B$ model perspective.

This same elementary structure is also present at genus 1:
\begin{conjecture}
The genus 1 Gromov-Witten potential of $X_k$ with the antidiagonal torus action is given in $B$ model variables by the closed formula
\ba
\no
G_k=\frac{11}{24}\log (1+(-1)^{k+1}(k+1)^2 q)+\Big(-\frac{5}{12}+\frac{(k+1)^2}{24}\Big)\log (1+(-1)^{k+1} q)-\frac{1}{2}\log \big(q\frac{d t_k}{dq}\big)
\ea
\end{conjecture}
This formula was derived, in a slightly different form, in \cite{M}, by
looking directly at the $A$ model calculation. Here, we have found this
through mirror symmetric methods, by interpreting the singular points of
the mirror map as the discriminant locus. 

These conjectures were arrived at by use of the following proposition:
\begin{proposition}
The equivariant Gromov-Witten invariants of 
\ba
X_k=\oo(k)\oplus\oo(-2-k)\rightarrow \p^1 
\ea
with action $(\la_1,\la_2)$ on the bundle are the same as those of the total space
\ba
X_k'=\oplus_1^k \oo(1)\bigoplus \oplus_1^{2+k}\oo(-1)\rightarrow \p^1
\ea
with $(\overbrace{\la_1\dots \la_1}^k;\overbrace{\la_2 \dots \la_2}^{2+k})$ acting on the bundle. 
\end{proposition}
This `factorization' of the bundle of $X_k$ into a sum of $\oo(1)$ and
$\oo(-1)$ terms has the effect of dramatically simplifying the mirror
map in the case of the antidiagonal action. We
will see that this proposition is actually the natural generalization of
the results found in \cite{FJ3}. Thus, the `mirror map' we have found
above for the antidiagonal torus action is actually the mirror map of
$X_k'$, rather than $X_k$.

To see why the above proposition is true, we
have the following argument due to H. Iritani. Let $p$ be the K\"ahler
class measuring the volume of the $\p^1$. Then we compare the
equivariant Euler
class of the bundle of $X_k$:
\ba
&e_T(\oo(k)\oplus \oo(-2-k))=(kp +\la_1)((-2-k)p+\la_2)=\\
&\big(k\la_2-(2+k)\la_1\big)p+\la_1\la_2
\ea
to that of $X_k'$:
\ba
&e_T(\oplus_1^k \oo(1)\bigoplus \oplus_1^{2+k}\oo(-1))=(p
+\la_1)^k(-p+\la_2)^{2+k}=\\
&\Big(\big(k\la_2-(2+k)\la_1\big)p+\la_1\la_2\Big)\la_1^{k-1}\la_2^{2+k-1}
\ea
where we have imposed the cohomology relation $p^2=0$.
Then we see that these are effectively the same (up to the multiplicative
factor $\la_1^{k-1}\la_2^{2+k-1}$, which can be thought of as simply a
product of trivial bundles). As the theorem of \cite{CG} uses only this
equivariant Euler class as input, one concludes that the equivariant
Gromov-Witten invariants of both spaces must be equal.  

Although we can prove the proposition this way, this equivalence was
originally deduced through more geometric considerations, which are described in
the body of the paper. Yet even with this proposition in hand,
it is quite nontrivial that we find the
remarkable structure listed above for the antidiagonal action,
i.e. the simplified mirror map and, most dramatically, the existence of
the rational Yukawa coupling.

The diagonal action $\la_1=\la_2$ unfortunately does not possess the nice structure observed above. This is not terribly surprising, as this case is manifestly not Calabi-Yau. Nonetheless, we find the following phenomenon:

\begin{conjecture}
The equivariant mirror symmetry computation on $X_k$ with the diagonal torus action is the same as that of $\p(\oo\oplus\oo(k)\oplus\oo(-2-k))$. That is, the mirror maps and Gromov-Witten invariants are equal. 
\end{conjecture}

We can gain some understanding of how this comes about by examining the charge vectors of the projective bundle spaces:
\ba
\p(\oo\oplus\oo(k)\oplus\oo(-2-k))&:&
\begin{pmatrix}
1&1&k&-2-k&0 \\
0&0&1&1&1
\end{pmatrix}
\ea
By looking at this matrix and then back at (\ref{matrix0}), we see that the toric data of these two spaces is essentially equivalent, and since the $I$ functions are determined entirely from the above matrices, it is not too outlandish that we should find the same mirror maps and Gromov-Witten invariants between these examples.

\subsection{$A_n$ singularities and toric trees}

In \cite{FJ3}, it was suggested that the Gromov-Witten invariants that physicists use are often actually equivariant Gromov-Witten invariants. Here, we realize this idea by computing the prepotential and genus 1 Gromov-Witten potential for $A_n$ singularities and for threefolds $X$ satisfying $\dim H_4(X,\Z)=0$. 

In \cite{FJ2}, we showed that one could compute the prepotential of a noncompact Calabi-Yau threefold up to polynomial terms of degree 2 by using various compactifications. However, particularly in the no 4 cycle $\dim H_4(X,\Z)=0$ case, this approach is not satisfactory, since we deliberately use a compactification known to reproduce the physically expected answer. Here, we will see that through the equivariant formalism, physical Gromov-Witten invariants appear naturally. 

We first consider the $A_n$ singularity. This geometry is described by the $n \times (n+2)$ matrix 
\ba
\begin{pmatrix}
1&-2&1&0&0&\cdots&0 \\
0&1&-2&1&0&\cdots&0\\
\vdots \\
0&0&0&\cdots&1&-2&1
\end{pmatrix}
\ea
This space has $n$ 2 cycles arranged along a line corresponding to the Dynkin diagram of $A_n$; we label these sequentially by $C_1\dots C_n$. Let $t_1\dots t_n$ be the complexified K\"ahler classes corresponding to these curves.
Then we work with the equivariant theory 
\ba
\begin{pmatrix}
0&\la_1&\la_2&\cdots&\la_{n-1}&\la_n&0\\
1&-2&1&0&0&\cdots&0 \\
0&1&-2&1&0&\cdots&0\\
\vdots \\
0&0&0&\cdots&1&-2&1
\end{pmatrix}
\ea
The idea behind this choice is the same as that of the $\oo\oplus \oo(-2)\rightarrow \p^1$ case considered in \cite{FJ3}: each $-2$ entry corresponds to a noncompact divisor, so we `compactify' these divisors by adding in an equivariant parameter for each. Then, we simply use the equivariant $I$ function to extract Gromov-Witten invariants via mirror symmetry. The result for the prepotential is the following:
\ba
\mathcal F_{A_n}=\sum_{i=1}^n Li_3(e^{t_i})+\sum_{i=1}^{n-1}Li_3(e^{t_i+t_{i+1}})+\dots+Li_3(e^{t_1+\dots+t_n}).
\ea
Later in the paper, we use this instanton expansion, together with the discriminant locus computed from a $\p^1$ fibration over the $A_2$ singularity, to exhibit mirror symmetry at genus 1 on $A_2$.

Now let $X$ be any noncompact Calabi-Yau threefold with no 4 cycles, such that $X$ is described by symplectic reduction via a matrix $M_{ab}$. As is well known, by making appropriate choices of curves $C_1\dots C_j$ and divisors $D_1\dots D_k$ in $X$, the entries of $M_{ab}$ give intersection numbers between curves and divisors. Suppose that $D_1\dots D_l$ is a basis of noncompact divisors of $X$. Then as above, we consider the equivariant Gromov-Witten theory of $X$ with exactly one equivariant parameter inserted for each noncompact divisor: 
\ba
\begin{pmatrix}
\la_1 \dots & \la_l& 0 \dots & 0 \\
m_{1,1}\dots & m_{1,l}& m_{1,l+1} \dots &m_{1,k} \\
\vdots \\
m_{j,1}\dots & m_{j,l}& m_{j,l+1} \dots &m_{j,k}
\end{pmatrix}
\ea
Again, only the equivariant $I$ function is required to compute Gromov-Witten invariants. The result for the prepotential is nearly the same as the above: we get one term in the prepotential for each curve in the geometry. However, there is one important difference: by tuning the equivariant parameters, we can arrange things so that the curves with normal bundle $\oo(-1)\oplus \oo(-1)$ and those with normal bundle $\oo \oplus \oo(-2)$ have either the same relative sign, or the opposite relative sign. Moreover, the choice in which the curves have opposite relative sign corresponds to the physical `anti-diagonal action' case, which is consistent with physically computed prepotentials. However, from the equivariant point of view, either sign convention is equally acceptable, as was suggested in \cite{FJ2}.

As the prepotential cannot be written in a concise form, we will instead
work out the explicit example of the trivalent $(-1,-1)$ curve in the text, computing the prepotential and subsequently exhibiting genus 1 mirror symmetry.

\section{Equivariant local mirror symmetry of curves}

\subsection{Review of previous results}

We begin with an overview of the findings of \cite{FJ3}. What was shown was essentially that we can use the equivariant version of the Givental $I$ function to compute equivariant Gromov-Witten invariants of $X_k$ with the general torus action $(\la_1,\la_2)$ on the bundle. For $X_{-1}$, the equivariant $I$ function reads 
\ba
I^{\la}_{-1}=e^{p\log q/\hbar}\sum_{d\ge 0}\frac{\prod_{i=1,2}\prod_{m=-d+1}^0(-p+m\hbar+\lambda_i)}{\prod_{m=1}^d(p+m\hbar)^2}q^d
\ea  
and on $X_0$,
\ba
I^{\la}_{0}=e^{p\log q/\hbar}\sum_{d\ge 0}\frac{\prod_{m=-2d+1}^0(-2p+m\hbar+\lambda)}{\prod_{m=1}^d(p+m\hbar)^2}q^d.
\ea 
These $I$ functions are annihilated by the following two equivariant differential operators, respectively:
\ba
\label{ops}
\mathcal D_{-1}&=&\theta^2-q(\theta-\la_1)(\theta-\la_2),\\
\mathcal D_0&=&\theta^2-q(2\theta-\la)(2\theta-\la+\hbar).
\ea
where $\theta=\hbar q d/dq$.
As was shown in \cite{FJ3}, these two $I$ functions agree up to the
mirror map and equivariant mirror map if we take $\la_1=\la_2=\la$ in
$I^{\la}_{-1}$. This means that the two equivariant differential
equations $\mathcal D_{-1}f=0,\mathcal D_0 f=0$ generate the same
quantum cohomology ring when $\la_1=\la_2$.

We note one unusual feature of $I^{\la}_0$ which will be key to the
derivations that follow. The issue is that the $I$ function is unable to
detect the trivial $\C$ factor of $X_0=\oo\oplus \oo(-2)\rightarrow
\p^1$. This implies two things. First, the $I$ function cannot be used
to compute the minus sign on the instanton expansion of $X_0$ claimed by
physicists for the antidiagonal action (e.g. \cite{M}). Secondly, the
point which central to this paper: $I^{\la}_0$ is actually the
equivariant $I$ function of the $A_1$ singularity. In other words, the
equivalence between $I^{\la_1=\la_2}_{-1}$ and $I^{\la}_0$ in \cite
{FJ3} is an equivalence of equivariant quantum cohomology rings between
the three dimensional space $X_0$ and the two dimensional $A_1$
singularity. The natural generalization of this observation is the Proposition of the previous section. We will see how this can be derived for $\oo(1)\oplus \oo(-3)\rightarrow \p^1$ with the antidiagonal action later in the paper. 

As one final remark on the above, we compare the equivariant charge vectors of $X_{-1}$ and the $A_1$ singularity:
\ba
\label{A_1}
X_{-1}&:&
\begin{pmatrix}
1&1&-1&-1 \\
0&0&\la &\la
\end{pmatrix}
\\
A_1&:&
\begin{pmatrix}
1&1&-2 \\
0&0&\la
\end{pmatrix}
\ea
From this vantage, the calculation on $X_{-1}$, while equivalent to that
of $A_1$, is slightly simpler, because there is no mirror map. So, the
dimension of the space has gone up, and the complexity of the mirror map
has gone down. This is the first example of the `factorization of the
bundle' stated in Proposition 1.

Moving on, from \cite{FJ3} we have the equivariant $I$ function for $X_k$ for $k\ge 1$:
\ba
I_k^{\la}=e^{p\log q/\hbar}\sum_{d\ge 0}\frac{\prod_{m=(-2-k)d+1}^0((-2-k)p+m\hbar+\lambda_2)}{\prod_{m=1}^d(p+m\hbar)^2 \prod_{m=1}^{kd}(kp+m\hbar+\lambda_1)}q^d.
\ea
Then we proposed that the equivariant Gromov-Witten invariants of $X_k$ could be recovered by first expanding $I_k^{\la}$ about $\la_1=\infty$, then performing Birkhoff factorization of the result (to remove positive powers of $\hbar$), and finally by inverting the mirror map and equivariant mirror map of the Birkhoff factorized function (which we called the $J$ function). 

We briefly recall the Birkhoff factorization procedure. Since the only
examples we work with will be curves, we can give an
especially simple formulation. Suppose we have an equivariant $I$ function
representing some bundle over a curve, and that after expansion about
the equivariant parameters $\la=\infty$, we obtain $I\in
\C[\hbar,\hbar^{-1}]$. Since such a power series expasion strictly
speaking does not make sense, we have to remove positive powers of
$\hbar$ from the $I$ function before extracting mirror symmetry
data. This is done by a theorem in \cite{CG}: there exist functions
$c_0,c_1$ such that
\ba
c_0(q,\hbar)I(q,\hbar,\hbar^{-1})+c_1(q,\hbar)\hbar q \frac{d}{dq}I(q,\hbar,\hbar^{-1})=J(q,\hbar^{-1})
\ea
and $J$ is independent of $\hbar$. We then obtain the mirror map,
etc. by looking at the $\hbar^{-1}$ expasion of $J$.

Then the above process of equivariant mirror symmetry, given in slightly
more detail, proceeds by performing the series expansion and Birkhoff
factorization, from which we find
\ba
J=1+\frac{t_1^{\la}(q) p+t_2^{\la}(q)}{\hbar}+\frac{W_1^{\la}(q) p+W_2^{\la}(q)}{\hbar^2}+\dots
\ea
Then by multiplying $J$ by $e^{-t_2^{\la}(q)}$ and inverting the mirror map $ t_1^{\la}(q)$, we can read off the instanton information from the $W_i(q).$ It was then shown that if $\la_1=\la_2$, the resulting functions are independent of $k$, and for $\la_1=-\la_2$, physical Gromov-Witten invariants for $X_k$ could be computed.

Unfortunately, this method produces no closed formulas, and requires serious computer power even to obtain results up to degree 6. Moreover, the resulting mirror maps have incredibly complicated formulas. For example, on the $X_1$ geometry with the diagonal torus action, the formula for the mirror map is given by
\ba
\frac{d \log q}{dt}=\frac{3}{8}\Big(1+\frac{\sin(\frac{5}{3}\sin^{-1}(\sqrt{108q}))}{\sqrt{108q}}\Big).
\ea 
These types of formulas have put the search for Yukawa couplings, as well as the genus 1 computation, out of reach. We were unable to even identify a discriminant for this case.  

These problems, and the desire to perform the $B$ model at genus 1, led us to search more carefully for the meaning of mirror symmetry for these spaces. The first clue is given by a close look at $D_n$ singularities.

\subsection{Connection between $X_1$ and $D$ type singularities}

In the previous section, we have explored the direct approach to mirror
symmetry for $X_k$ via $I$ functions, and have noted along the way that
the $X_0$ calculation can be viewed as equivariant theory on the $A_1$
singularity. Next, we claim that the equivariant Gromov-Witten
invariants of $X_1=\oo(1)\oplus \oo(-3)\rightarrow \p^1$ are the same as the
equivariant Gromov-Witten invariants of a certain partial resolution of
the $\C^2/D_n$ singularity.

Recall that a simple singularity $\C^2/D_n$,$n\ge 4$,where $D_n$ is the $n$th dihedral group, can be realized as a hypersurface in $\C^3$ \cite{G4}
\ba
f=x_1^2 x_2-x_2^{n-1}+x_3^2=0.
\ea
One can obtain a smooth variety in two ways. One way is complex deformation: take a basis $\{\mu_1\dots \mu_m\}$ of the local algebra of the singular point
\ba
H=\frac{\C[x_1,x_2,x_3]}{\langle \partial_{x_1}f,\partial_{x_2}f,\partial_{x_3}f\rangle}
\ea
and deform $f$ as $f^{\la}=f+\sum_{j=1}^m \la_j \mu_j$.
The second way is by blowing up the $n$ singular points; we end up with
$n$ curves of self-intersection $-2$, and moreover there is a special
 `central' curve which intersects $3$ other curves exactly once. We call
this the trivalent curve. Mirror symmetry for the $D_n$ singularity is
then realized as the transformation between the blown up space and the
complex deformed space.

We next recall the work of Cachazo-Katz-Vafa \cite{CKV}, where it was
shown (in the so-called Laufer's example section, p.37-40) that a
certain monodromic fibration of the $D_n$ singularity over the plane,
\it where only the trivalent curve is blown up, \normalfont is equivalent
to the geometry $X_1=\oo(1)\oplus \oo(-3)\rightarrow \p^1$. Now, in our
case, we can't use this fact directly, because we are interested in the
mirror symmetry computation. From the vantage of mirror symmetry,
everything is much simpler if we stay in the realm of toric
geometry. Nonetheless, there is reason to suspect that we can still find a relationship between $X_1$ and $D_n$ singularities at the toric level; see
for example \cite{K}, where the connection is described as follows. If
we let $X$ be a Calabi-Yau threefold with imbedded curve
$\p^1\hookrightarrow X$ such that $N_{\p^1/X}\cong \oo(1)\oplus
\oo(-3)$, and we then shrink the $\p^1$ to a point $p$, then under
certain conditions the
singularity type of the generic hyperplane section through the point $p$ will be of type $D_4$.

Thus, we would like to consider the toric representation of the $D_n$ geometry in which only the trivalent curve is blown up.
From \cite{KMV}, the Mori cone vector corresponding to this trivalent curve is given by
\ba
\label{v2}
\begin{pmatrix}
1&1&1&-2&-1
\end{pmatrix},
\ea
and therefore we end up with a geometry defined by single K\"ahler parameter.
We call this the `$D_1$ singularity'.

We now show that this is indeed the Mori cone generator of the trivalent
curve. Consider for simplicity the blown up $D_4$
geometry, and let $C_1\dots C_4$ be a basis of curve classes, where
$C_4$ is the trivalent curve. We claim that the toric data defining the
blowup of the
$D_4$ singularity is given by the matrix
\ba
\begin{pmatrix}
l_1\\l_2\\l_3\\l_4
\end{pmatrix}=
\begin{pmatrix}
-2&0&0&1&0&1&0&0\\
0&-2&0&1&0&0&1&0\\
0&0&-2&1&0&0&0&1\\
1&1&1&-2&-1&0&0&0
\end{pmatrix}
\ea

We can derive this matrix as follows. The rows of the matrix, as well as the first 4
columns, correspond to the curves
$C_1\dots C_4$, and entries of the matrix give intersection numbers
between curves and divisors in the geometry.
Hence the $-2$ entries are interpreted as the self-intersection numbers
$C_i^2=-2$, and e.g. the meaning of the $(4,1)$ entry is $C_1 \cdot
C_4=1$, which is true because $C_4$ intersects $C_1\dots C_3$ exactly once. Thus we
see the necessity of the three 1's in the fourth row. Also, the last 3
columns represent the (noncompact) normal bundles to the curves
$C_1\dots C_3$ respectively, and since the intersection number of each
curve with its normal bundle is $+1$, we obtain the entries of these columns. Finally, in order to
impose the Calabi-Yau condition on the space, we need to add the -1 in the fifth column (so
$\sum_j l_4^j=0$). Hence we arrive at the claimed form of the Mori
vector (\ref{v2})
  
Since this is the matrix of intersection numbers of the space, we can
represent blown up $D_4$ as a complex 4 dimensional space given by 
\ba
\{(z_1\dots z_8)\in \C^8:\sum_{j=1}^8 l_i^j |z_j|^2=r_i, i=1\dots 4\}/(S_1)^4
\ea
where $r_i$ are real parameters and the action is given by 
\ba
S^1_i:(z_1\dots z_8)\longrightarrow(e^{l_i^1 \sqrt{-1}\theta_1}z_1\dots e^{l_i^8 \sqrt{-1}\theta_8}z_8)
\ea

We now return to the discussion of the relation between the $D_1$
geometry and that of $X_1$. Notice that the charge vector of $X_1=\oo(1)\oplus
\oo(-3)\rightarrow \p^1$ is given as:
\ba
\label{v1}
\begin{pmatrix}
1&1&1&-3
\end{pmatrix}.
\ea
If we stare at the vectors in Eqns (\ref{v1}), (\ref{v2}) and then look back at Eqn.(\ref{A_1}), it is not hard to imagine that the following two equivariant theories would give the same Gromov-Witten invariants:
\ba
D_1&:&
\begin{pmatrix}
1&1&1&-1&-2 \\
0&0&\la_1&\la_2&\la_2
\end{pmatrix}
\\
X_1&:&
\begin{pmatrix}
1&1&1&-3 \\
0&0&\la_1&\la_2
\end{pmatrix}
\ea
Indeed, direct computation verifies the equality of the Gromov-Witten
invariants. We note that as in the $X_{-1}$/ $A_1$ case, the dimension
of the spaces is different. The charge vector of the $D_1$ geometry
identifies the space as $\oo(1)\oplus \oo(-1)\oplus \oo(-2)\rightarrow
\p^1$, a complex fourfold, which is consistent with the fourfold
representation obtained for the $D_4$ singularity above. 

However, the mirror map on $D_1$ and that on $X_1$ are not the
same. What does this mean? Can we hope that the mirror map is somehow
getting simpler with increased dimension, as on the $X_{-1}\cong A_1$ example?

The answer, as well as the means of producing rational Yukawa couplings, lies in specializing the torus weights to the equivariantly Calabi-Yau setting $\la_1=-\la_2=\la$. Now consider the differential operator which annihilates the $I$ function of the $D_1$ singularity:
\ba
\mathcal D_{D_1}=\theta^2(\theta+\la)-q(-\theta-\la)(-2\theta-\la)(-2\theta-\la-\hbar)
\ea
This operator can be factorized as
\ba
\mathcal D_{D_1}=\big(\theta^2+q(-2\theta-\la)(-2\theta-\la-\hbar)\big)(\theta+\la)
\ea
Then Eqn.(\ref{ops}) asserts the equivalence of the operators
$\theta^2+q(-2\theta-\la)(-2\theta-\la-\hbar)$
and $\theta^2+q(\theta+\la)^2$, from which we expect that we can also reproduce equivalent Gromov-Witten invariants by use of the operator
\ba
\mathcal D'&=&(\theta^2+q(\theta+\la)^2)(\theta+\la)\\
&=&\theta^2(\theta+\la)-q(-\theta-\la)^3.
\ea
This last form of the $\mathcal D'$ operator can be derived from the toric data
\ba
\begin{pmatrix}
1&1&1&-1&-1&-1 \\
0&0&\la&-\la&-\la&-\la
\end{pmatrix}
\ea
and corresponds to the total space $\oo(1)\oplus\oo(-1)\oplus\oo(-1)\oplus\oo(-1)\rightarrow \p^1$ with a torus action $(\la,-\la,-\la,-\la)$ on the bundle. Again, we may directly compute to verify that indeed the equivariant Gromov-Witten invariants corresponding to this toric data agree with those of $X_1$ with the antidiagonal action.

Now we ask what we have gained through these geometric manipulations. We consider the $I$ function which generates the solution space of $\mathcal D' f=0$:
\ba 
I'=e^{p \log q / \hbar}\sum_{d\ge 0}\frac{\prod_{i=-d+1}^0 (-p-\la+m\hbar)^3}{\prod_{i=1}^d(p+\la+m\hbar)\prod_{i=1}^d(p+m\hbar)^2}q^d
\ea
We run $I'$ through the same procedure described in Section 3.1: expand $I'$ about $\la=\infty$ and then perform Birkhoff factorization of the result. After this, we find the following incredible result:
\ba
\label{easyJ}
J'=1+\frac{p\big(\log q+ 3\log (1+q)\big)+\la \log(1+q)}{\hbar}+\dots
\ea
In other words, the mirror map has taken on the nearly trivial form
$t=\log q+ 3\log (1+q)$, and moreover the equivariant mirror map (the
coefficient of $\la$ in Eqn.(\ref{easyJ})) is just a multiple of the
regular mirror map.

\subsection{Natural  rational $B$ model Yukawa couplings}

We now argue that the natural choice of the triple intersection number
of the curve, from the perspective of the $B$ model calculation, should
be $\frac{1}{k(k+2)}$. This will give the simple form of the
Picard-Fuchs equation and $B$ model Yukawa coupling of Conjecture 1. 
Moreover, this agrees with the predictions of \cite{AOSV}.

We consider the prepotential of $X_k$ with antidiagonal action and
arbitrary triple intersection number $\langle C,C,C \rangle$. Using the expansion of \cite{M}, this is
\ba
\mathcal F_k=\langle C,C,C \rangle\frac{t^3}{3!}-\sum_{d=1}^{\infty}\frac{(-1)^{kd}((k+1)^2d-1)!}{d!d^2(((k+1)^2-1)d)!}e^{td}.
\ea
The key to choosing the right intersection number now actually lies in
looking at the second derivative of this function. An easy computation
shows that
\ba
\frac{\partial^2 \mathcal F_k}{\partial t_k^2}=\langle C,C,C \rangle \log \big(q(1+(-1)^{k+1}
q)^{k(k+2)}\big)+\log \big(1+(-1)^{k+1}q\big).
\ea
Then by setting the value
\ba
\langle C,C,C \rangle=\frac{-1}{k(k+2)},
\ea
we immediately find the remarkably simple relation 
\ba
\frac{\partial^2 \mathcal F_k}{\partial t_k^2}=\frac{-1}{k(k+2)}\log q.
\ea

Rather than being merely the choice of intersection number which gives
the simplest expression of the above function, there are several
other reasons to expect that this is the correct choice (besides the
fact that it was predicted to be so in \cite{AOSV}). The first is by
examining the form of the Picard-Fuchs equation which results from the
above Yukawa coupling. 
The Picard-Fuchs equation describing mirror symmetry will be given
by 
\ba
 \partial_{t_k}\Big(\frac{1}{Y_{t_k t_k t_k}}\Big)\partial^2_{t_k}&=&\theta
k(k+2)q\frac{dt_k}{dq}\Big(q\frac{dt_k}{dq}\Big)^{-1}\theta
\Big(q\frac{dt_k}{dq}\Big)^{-1} \theta 
\\ &=& \theta^2
\Big(q\frac{dt_k}{dq}\Big)^{-1} \theta.
\ea
In other words, only the information of the mirror map is required to
compute all Gromov-Witten invariants of $X_k$. Said differently, mirror
symmetry is completely characterized by the integer $k$ (since the
mirror map is). 

Another reason the above choice of intersection number is the most
natural can be seen from the genus 1 function in the next section.
First, note that by using the above expression for $\partial^2
\mathcal F / \partial t^2$, we have the $B$ model Yukawa
coupling, as given in Conjecture 1:
\ba
Y_{qqq}=\frac{\partial^3
\mathcal F}{\partial t^3} \Big(q \frac{d t_k}{dq}\Big)^3= \frac{-1}{k(k+2)}\Big(q\frac{dt_k}{dq}\Big)^2.
\ea

As is
evident from the expansions of the next section, there are two
components of the discriminant locus, given by the numerator and
denomenator of the mirror map, respectively. Now, if we make any other
choice of triple intersection number, what happens is that the $B$ model
Yukawa coupling will contain an extra polynomial factor. From the
perspective of period integrals, this suggests that there will be an
extra component of the discriminant locus, which one would reasonably
expect to appear in the genus 1 function. However, as this extra
component does not appear, it is natural that the Yukawa coupling should
have the same singular points as the mirror map, thus lending support to
our choice of intersection number.

Finally, we can also use the formula proposed in \cite{FJ2} for the
computation of triple intersection numbers. Adapted to the present case,
this reads
\ba
\langle C,C,C \rangle =\int_{\p^1} \frac{J^3}{e(\oo(k)\oplus \oo(-2-k))}=\int_{\p^1}\frac{J^3}{J^2k(-2-k)}=-\frac{1}{k(2+k)}\int_{\p^1}J
\ea
where $J$ is the K\"ahler class satisfying $\int_{\p^1}J=1$ and $e()$
denotes the Euler class. Hence we
obtain the intersection number claimed.

\subsection{Genus one.}

With such a simple form for the Yukawa coupling and mirror map at hand, it is natural to suppose that we also have an elementary form for the genus 1 Gromov-Witten potential of $X_k$ on the $B$ model, which we denote by $G_k$. In fact, our interest in doing the genus 1 computation was the original motivation behind this project.

 The $A$ model function can be worked out by the topological vertex
 \cite{AKMV}. We are grateful to Y. Konishi for providing us with a
 program for the vertex calculation. Then, we need only compare the general form of the $B$ model genus 1 amplitude to see if it agrees with the $A$ model answer. Recall \cite{BCOV} that the $B$ model function $G_k$ has the general structure
\ba
G_k=\log \big(q^a\prod_i \Delta_i^{b_i} J_k \big)
\ea
where $\Delta_i$ is and irreducible component of the discriminant locus, $J_k=d \log q/dt_k$ is the Jacobian and $a,b$ are rational numbers. We again specialize to the case $X_1$ for clarity. Then we have $J_1=(1+q)/(1+4q)$, and since this is the derivative of the mirror map, the singular points of $J_1$ define the discriminant locus: 
\ba
\Delta_1=1+q, \ \ \ \ \ \ \Delta_2=1+4q.
\ea
Then a simple calculation verifies that $G_1$ is given as
\ba
G_1(q)=\frac{11}{24}\log (1+4q)-\frac{1}{4}\log(1+q)+\frac{1}{2}\log (J_1).
\ea
By substituting the inverse mirror map into $G_1$, we find the expansion 
\ba
G_1(t)=\frac{1}{12}e^t-\frac{1}{24}e^{2t}-\frac{29}{36}e^{3t}+\frac{499}{48}e^4t-\frac{517}{5}e^{5t}+\dots
\ea
In exactly the same way, we find $G_2$:
\ba
G_2=\frac{11}{24}\log (1-9q)-\frac{1}{24}\log(1-q)+\frac{1}{2}\log (J_2)
\ea
where now $J_2=(1-q)/(1-9q)$.
This has expansion 
\ba
G_2(t)=-\frac{1}{12}e^t+\frac{19}{24}e^{2t}+\frac{899}{36}e^{3t}+\frac{27259}{48}e^4t+\frac{733289}{60}e^{5t}+\dots
\ea

 By checking a few more cases we have formulated Conjecture 2.

\subsection{A word about the diagonal action.}

We now briefly consider equivariant mirror symmetry for $X_k$ with the
diagonal torus action $\la_1=\la_2$. We will see that this is very
likely the same calculation as the projective bundle $Y_k=\p(\oo\oplus
\oo(k)\oplus \oo(-2-k))\rightarrow \p^1$. Unfortunately, the techniques
of the earlier sections don't give any simple form for the mirror map or
genus 1 expansion, but as $Y_k$ is not Calabi-Yau, this is not entirely unexpected.

We see some hint of the correspondence already between the spaces
$X_{-1}$ and $Y_{-1}$. We have the expansion for the equivariant $I$
function on $X_{-1}$:
\ba
I_{X_{-1}}&=&e^{p\log q/\hbar}\sum_{d\ge
0}\frac{\prod_{m=-d+1}^0(-p+m\hbar+\la)}{\prod_{m=1}^d(p+m\hbar)}q^d\\
&=&1+\frac{p\log q}{\hbar}+\frac{\la^2 Li_2(q)-2p\la Li_2(q)}{\hbar^2}+\dots
\ea 
and that on $Y_{-1}$:
\ba
\no
I_{Y_{-1}}=e^{(p_1 \log q_1+p_2 \log q_2)/\hbar}\sum_{d\ge
0}\frac{\prod_{m=-\infty}^0(-p_1+p_2+m\hbar)}{\prod_{m=1}^{d_1}(p_1+m\hbar)^2\prod_{m=-\infty}^{-d_1+d_2}(-p_1+p_2+m\hbar)\prod_{m=1}^{d_2}(p_2+m\hbar)}q_1^{d_1}q_2^{d_2}\\ \no
=1+\frac{p_1 \log q_1+p_2 \log
q_2}{\hbar}+\frac{p_2^2Li_2(q_1)-2p_1p_2Li_2(q_1)+\log q_1 \log q_2+(\log q_2)^2/2}{\hbar^2}+\dots
\ea
The expansions disagree at higher order in $1/\hbar$, but one might
simply attribute this to the fact that $I_{Y_{-1}}$ is considered as a
cohomology valued hypergeometric series, taking values in the cohomology
ring
\ba
\frac{\C[p_1,p_2]}{\langle p_1^2, p_2(p_1-p_2)^2 \rangle},
\ea
while the coefficients of $I_{X_{-1}}$ are only subjected to the
relation $p^2=0$. This argument is strengthened if one looks at the next
example:
\ba
I_{X_0}=1+\frac{p(\log q +2f(q))-\la f(q)}{\hbar}+\dots,
\ea
\ba
I_{Y_0}=1+\frac{p_1(\log q_1 +2f(q_1))+p_2(\log q_2- f(q_1))}{\hbar}+\dots,
\ea
where $f(q)=\sum_{n>0}q^n(2n-1)!/(n!)^2$. That is, the mirror maps are
the same, and moreover the instanton expansions agree exactly after
inversion of the respective mirror maps. 

The more nontrivial statement is that this correspondence holds even
across Birkhoff factorization. Namely, if we consider $X_1$ with the
diagonal torus action $(\la,\la)$ and carry out Birkhoff factorization
on the resulting $I$ function, both the mirror map and the Gromov-Witten
invariants turn out to be exactly the same as those of the projective
bundle $Y_1=\p(\oo\oplus \oo(1)\oplus \oo(-3))\rightarrow \p^1$, where we
have performed Birkhoff factorization on the $I$ function for $Y_1$.
 
Unfortunately, we were unable to obtain nice formulas for any of the
quantities discussed in this paper on $Y_k$. Hence, the analog of
Picard-Fuchs equations, etc. remains unclear.

\section{Equivariant mirror symmetry for $X$ with $\dim H_4(X)=0$}

We now show that equivariant techniques can be used to effectively
compute Gromov-Witten potentials not only for bundles over curves, but
in fact for any Calabi-Yau lacking four cycles.

\subsection{$A_2$}

Consider the toric charge vectors for the standard $A_2$ geometry
\ba
\begin{pmatrix}
1&-2&1&0\\
0&1&-2&1
\end{pmatrix}
\ea
This has two curve classes $C_1,C_2$ corresponding to the rows of the
above matrix. Also let $p_1,p_2$ be K\"ahler classes satisfying $\int_{C_i}p_j=\delta_{ij}$. Note that for $C_1$, the normal bundle direction is given
by the second column of the matrix, and for $C_2$, the normal bundle
corresponds to the third column. As such, $\la_1$ is the equivariant
parameter corresponding to the $C_1$ curve, and $\la_1$ to $C_2$.

Now, in order to exhibit equivariant Gromov-Witten invariants on the
$A_1$ singularity, we simply added an equivariant parameter
corresponding to the normal bundle direction:
\ba
\begin{pmatrix}
1&1&-2
\end{pmatrix}\rightarrow
\begin{pmatrix}
0&0&\la \\
1&1&-2
\end{pmatrix}.
\ea
Hence the most natural equivariant theory we can use to extract
Gromov-Witten invariants is   
\ba
\label{matrix1}
\begin{pmatrix}
0&\la_1&\la_2&0\\
1&-2&1&0\\
0&1&-2&1
\end{pmatrix}
\ea
The strategy is then simple enough: just use the equivariant $I$
function to compute the mirror map and prepotential. The $I$ function is 
\ba
I_{A_2}=e^{(p_1 \log q_1 +p_2 \log q_2)/\hbar}\sum_{d \ge 0}C(d_1,d_2,\la)q_1^{d_1}q_2^{d_2}
\ea
where $C(d_1,d_2,\la)=$
\ba
\no
\frac{\prod_{m=-\infty}^0(-2p_1+p_2+\la_1+m\hbar)\prod_{m=-\infty}^0(p_1-2p_2+\la_2+m\hbar)}{\prod_{m=1}^{d_1}(p_1+m\hbar)\prod_{m=-\infty}^{-2d_1+d_2}(-2p_1+p_2+\la_1+m\hbar)\prod_{m=-\infty}^{d_1-2d_2}(p_1-2p_2+\la_2+m\hbar)\prod_{m=1}^{d_2}(p_2+m\hbar)}
\ea
Recall that these coefficients are subjected to the cohomology relations $p_1^2=p_1p_2=p_2^2=0$.
Then we go through the usual motions of expanding this function in
powers of $1/\hbar$ and inverting the mirror map and equivariant mirror
map, which are
given by the coefficient of $1/\hbar$. Let $t_1,t_2$ be the mirror map,
and set $x_1=e^{t_1}$,$x_2=e^{t_2}$.
Let $J_{A_n}$ be the function obtained by coordinate change of $I_{A_n}$ by
the mirror map. Then we want to read the instanton information from the
coefficient $W$ of $1/\hbar^2$ of $J_{A_n}$. There is one minor point one
needs to keep in mind when extracting the instanton information: since
the $\la_1$ equivariant parameter corresponds to the curve $C_1$, there
are certain `anomalous' terms in $W$ which mix the normal bundle of $C_1$ and the curve $C_2$. We
can cancel these terms by looking at one equivariant parameter at a
time: if we set $p_2=\la_2=0$, we obtain
\ba
\label{supera}
W|_{\la_2=0,p_2=0}&=&\la_1^2\Big(Li_2(x_1)+Li_2(x_1x_2)\Big)+p_1\la_1\Big(2Li_2(x_1)+Li_2(x_1x_2)-Li_2(x_2)\Big)\\
&=&\la_1^2\frac{\partial \mathcal F}{\partial
t_1}+p_1\la_1\Big(2\frac{\partial \mathcal F}{\partial
t_1}-\frac{\partial \mathcal F}{\partial
t_2}\Big)
\ea
where
\ba
\mathcal F=Li_3(x_1)+Li_3(x_2)+Li_3(x_1x_2).
\ea
Note that this is just as one would expect from local mirror symmetry
calculations, since the coefficient of $p_1\la_1$ of Eqn.(\ref{supera})
features a linear combination of prepotential derivatives determined
by the second column of the matrix (\ref{matrix1}).

Now that we know the instanton expansion, we can use this to work out
mirror symmetry at genus 1 for the $A_2$ singularity. Unfortunately, the
geometry of the mirror of (\ref{matrix1}) is degenerate, which means we
can't easily extract the discriminant locus from the mirror
manifold. Hence, we instead use a $\p^1$ fibration over $A_2$
\ba
\begin{pmatrix}
1&1&-2&0&0&0\\
0&0&1&-2&1&0\\
0&0&0&1&-2&1
\end{pmatrix}
\ea
We then compute the discriminant locus from the mirror manifold, and
then take the limit as the first curve disappears. The result is
\begin{eqnarray}
\no
\Delta=&&1-8q_1-8q_2+68q_1q_2+16q_1^2+16q_2^2-144q_1q_2^2-144q_1^2q_2+270q_1^2q_2^2+216q_1^3q_2^2+216q_1^2q_2^3\\
\no
&&-972q_1^3q_2^3+729q_1^4q_2^4,
\end{eqnarray}
from which we can immediately exhibit genus 1 mirror symmetry:
\ba
G_{A_2}=\frac{t_1}{12}+\frac{t_2}{12}-\frac{1}{12}\log
\big((1-e^{t_1})(1-e^{t_2})(1-e^{t_1+t_2})
\big)=\log\big(q_1^{1/12}q_2^{1/12}\Delta^{-7/24}\big(\frac{\partial \log
q}{\partial t}\big)^{1/2}\big).
\ea 
$\partial \log q/\partial t$ is the Jacobian of the mirror map.
The coefficient of $-7/24$ is also the same for the $A_1$ singularity,
suggesting that this behavior may be universal for $A_n$ singularities.

\subsection{The trivalent $(-1,-1)$ curve}

Since this example is closely related to the above, we give only the
briefest discussion, merely indicating the points at which this differs
from $A_n$. The geometry we are considering has 3 curves with
normal bundle $\oo(-1)\oplus \oo(-1)$ and all 3 curves intersect at a
single point. The equivariant theory we use is thus
\begin{eqnarray}
\begin{pmatrix}
0&0&0&\la_1&\la_2&\la_3\\
1&0&0&1&-1&-1\\
0&1&0&-1&1&-1\\
0&0&1&-1&-1&1
\end{pmatrix}
\end{eqnarray}
since the last 3 columns correspond to noncompact divisors. Then, as
above, we can use the equivariant $I$ function to work out mirror
symmetry. Let $W$ be the coefficient of $1/\hbar^2$ of the mirror map
transformed $J$ function. Then we restrict to the curve $C_1$
corresponding to the second row of the above matrix by setting
$\la_1=p_2=p_3=0$. If we choose the diagonal action, for which $\la_2=\la_3=\la$,
\ba
W|_{\la_1=p_2=p_3=0,\la_2=\la_3=\la}=\la^2\frac{\partial \mathcal
F_1}{\partial t_1}+p_1\la\big(2\frac{\partial \mathcal
F_1}{\partial t_1}-\frac{\partial \mathcal
F_1}{\partial t_2}-\frac{\partial \mathcal
F_1}{\partial t_3}\big)
\ea
and for the antidiagonal action, 
\ba
W|_{\la_1=p_2=p_3=0,\la_2=-\la_3=\la}=\la^2\frac{\partial \mathcal
F_2}{\partial t_1}+p_1\la\big(\frac{\partial \mathcal
F_2}{\partial t_3}-\frac{\partial \mathcal
F_2}{\partial t_2}\big)
\ea
where
\ba
\mathcal
F_k=\sum_{i=1}^3Li_3(x_i)+Li_3(x_1x_2x_3)+(-1)^{k+1}\big(Li_3(x_1x_2)+Li_3(x_1x_3)+Li_3(x_2x_3)
\big).
\ea
In other words, the use of the diagonal versus antidiagonal action
changes the relative sign between $(-1,-1)$ curves and $(0,-2)$
curves. Naturally, this means that neither sign choice is preferred
equivariantly, and indeed this continues to hold true at genus 1. 

\section{Conclusion}

In this paper, we have uncovered a surprisingly simple structure
underlying  mirror symmetry in genus 0 and 1 on 
$X_k=\oo(k)\oplus \oo(-2-k)\rightarrow \p^1$ with antidiagonal
action. Although one would expect some simplification of the calculation
using the factorization of the bundle as described in Proposition 1, the
miraculous appearance of a rational Yukawa coupling points toward some
deeper structure behind the problem. 

There are several possible directions for future work. One obvious
problem is the extension of these results to higher genus on the $B$
model, which involves the computation of the holomorphic anomaly at each
genus. One might also consider whether there is some nice form for the $B$ model
computation on $\p(\oo \oplus \oo(k)\oplus \oo(-2-k))$, or equivalently
$X_k$ with the diagonal action. We were able to derive rational $B$
model couplings for $\p(\oo \oplus \oo \oplus \oo(-2))$, but the
complexity of the mirror map for $k \ge 1$ poses a major
obstacle. Finally, one could also include open strings into the
computation. We expect to address these issues in future work.

\end{document}